\documentclass[letterpaper,final,12pt,reqno]{amsart}

\usepackage[total={6.3in,9.2in},top=1.1in,left=1.1in]{geometry}

\usepackage{times,bm,bbm,empheq,verbatim,fancyvrb,graphicx}
\usepackage[dvipsnames]{xcolor}

\usepackage[kw]{pseudo}

\pseudoset{left-margin=15mm,topsep=5mm,idfont=\texttt}

\usepackage{tikz}
\usetikzlibrary{decorations.pathreplacing}

\usepackage[pdftex,
colorlinks=true,
plainpages=false, 
linkcolor=blue,   
citecolor=Red,    
urlcolor=black    
]{hyperref}

\DefineVerbatimEnvironment{cline}{Verbatim}{fontsize=\small,xleftmargin=5mm}

\newcommand{\eps}{\epsilon}
\newcommand{\RR}{\mathbb{R}}

\newcommand{\grad}{\nabla}

\newcommand{\be}{\mathbf{e}}

\newcommand{\br}{\mathbf{r}}
\newcommand{\bu}{\mathbf{u}}

\newcommand{\ip}[2]{\left<#1,#2\right>}

\newcommand{\Rpr}{R_{\text{pr}}}
\newcommand{\Rin}{R_{\text{in}}}
\newcommand{\Rfw}{R_{\text{fw}}}

\begin{document}
\title[The FAS multigrid scheme]{The full approximation storage multigrid scheme: \\ A 1D finite element example}

\author{Ed Bueler}

\begin{abstract}  This note describes the full approximation storage (FAS) multigrid scheme for an easy one-dimensional nonlinear boundary value problem discretized by a simple finite element (FE) scheme.  We apply both FAS V-cycles and F-cycles, with a nonlinear Gauss-Seidel smoother, to solve the finite-dimensional problem.  The mathematics of the FAS restriction and prolongation operators, in the FE case, are explained.  A self-contained Python program implements the scheme.  Optimal performance, i.e.~work proportional to the number of unknowns, is demonstrated for both kinds of cycles, including convergence nearly to discretization error in a single F-cycle.  \end{abstract}

\thanks{Version 3.  This note is expository, and submission for publication is not foreseen.  Thanks to Matt Knepley for thoughtful comments.}

\maketitle

\tableofcontents

\thispagestyle{empty}

\section{Introduction}  \label{sec:intro}

We consider the full approximation storage (FAS) scheme, originally described by Brandt \cite{Brandt1977}, for an easy nonlinear elliptic equation.  Like other multigrid schemes it exhibits optimal solver complexity \cite{Bueler2021} when correctly applied, as we demonstrate at the end.  Helpful write-ups of FAS can be found in well-known textbooks \cite{BrandtLivne2011,Briggsetal2000,Trottenbergetal2001}, but we describe the scheme from a finite element point of view, compatible with the multigrid approaches used for obstacle problems \cite{GraeserKornhuber2009} for example, and we provide an easy-to-digest Python implementation.

Our problem is an ordinary differential equation (ODE) boundary value problem, the nonlinear Liouville-Bratu equation \cite{Bratu1914,Liouville1853}:
\begin{equation}
  -u'' - \lambda\, e^u = g,  \qquad u(0) = u(1) = 0.  \label{liouvillebratu}
\end{equation}
In this problem $\lambda$ is a real constant, $g(x)$ is given, and we seek $u(x)$.  This equation arises in the theory of combustion \cite{FrankKameneckij1955} and the stability of stars.

Our goal is to solve \eqref{liouvillebratu} in optimal $O(m)$ time on a mesh of $m$ elements.  A Python implementation of FAS, \texttt{fas1.py} in directory \texttt{fas/py/},\footnote{Clone the Git repository\, \href{https://github.com/bueler/fas-intro}{\texttt{github.com/bueler/fas-intro}}\, and look in the \texttt{fas/py/} directory.} accomplishes such optimal-time solutions both by V-cycle and F-cycle strategies (section \ref{sec:cycles}), and this note serves as its documentation.  While optimal-time solutions of 1D problems are not unusual, FAS and other multigrid strategies for many nonlinear 2D and 3D partial differential equations (PDEs) are also optimal.  This makes them the highest-performing class of solver algorithms for such problems.

By default the program \texttt{fas1.py} solves \eqref{liouvillebratu} with $g=0$.  A runtime option \texttt{-mms}, the ``method of manufactured solutions'' \cite{Bueler2021}, facilitates testing by specifying a problem with known exact solution and nonzero source term.  In detail, the solution is $u(x)=\sin(3\pi x)$, and by differentiation $g(x)=9\pi^2 \sin(3\pi x) - \lambda e^{\sin(3\pi x)}$.

\section{The finite element method}  \label{sec:femethod}

To solve the problem using the finite element (FE) method \cite{Braess2007,Bueler2021,Elmanetal2014}, we rewrite \eqref{liouvillebratu} in weak form.  Let $F$ be the nonlinear operator
\begin{equation}
  F(u)[v] = \int_0^1 u'(x) v'(x) - \lambda e^{u(x)} v(x)\, dx,  \label{operator}
\end{equation}
acting on $u$ and $v$ from the space of functions $\mathcal{H}=H_0^1[0,1]$, a Sobolev space \cite{Evans2010}.  (These functions have zero boundary values and one square-integrable derivative.)  Note $F(u)[v]$ is linear in $v$ but not in $u$.  We also define a linear functional built from the right-hand function $g$ in \eqref{liouvillebratu}:
\begin{equation}
  \ell[v] = \ip{g}{v} = \int_0^1 g(x) v(x) dx.  \label{rhsfunctional}
\end{equation}
Both $F(u)[\cdot]$ and $\ell[\cdot]$ are (continuous) linear functionals, acting on functions $v$ in $\mathcal{H}$, thus they are in the dual space $\mathcal{H}'$.  One derives the weak form
\begin{equation}
  F(u)[v] = \ell[v] \qquad \text{for all $v$ in $\mathcal{H}$} \label{weakform}
\end{equation}
by multiplying equation \eqref{liouvillebratu} by a test function $v$ and integrating by parts.  From now on we address problem \eqref{weakform}, despite its abstract form.

In an FE context a clear separation is desirable between functions, like the solution $u$, and the equations themselves, which are, essentially, functionals.  As in linear algebra, where one indexes the equations by row indices, \eqref{weakform} states the ``$v$th equation''; the equations are indexed by the test functions.  The FE method will reduce the problem to a finite number of unknowns by writing $u$ in a basis of a finite-dimensional subspace of $\mathcal{H}$.  One gets finitely-many equations by using test functions $v$ from the same basis.

\begin{figure}
\includegraphics[width=0.6\textwidth]{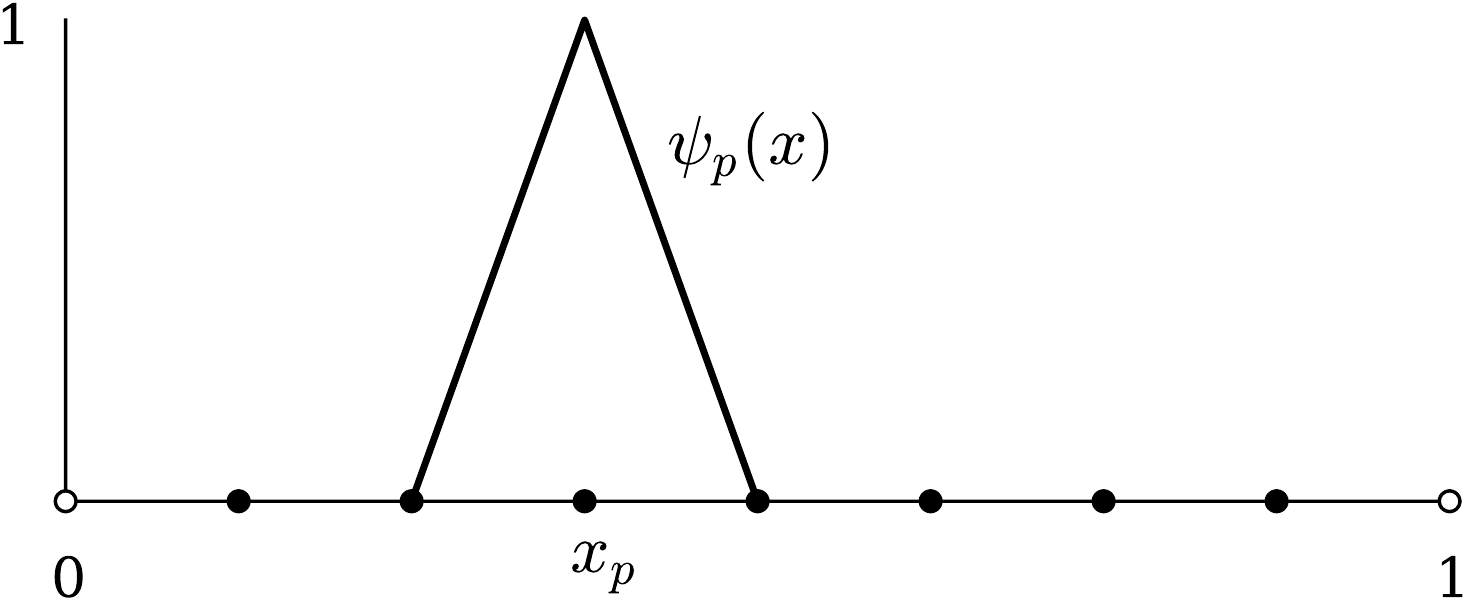}
\caption{A piecewise-linear hat function $\psi_p(x)$ lives at each interior node $x_p$.}
\label{fig:onehat}
\end{figure}

We apply the simplest possible mesh setup, namely an equally-spaced mesh on $[0,1]$ of $m$ elements (subintervals) of lengths $h=1/m$.  The interior nodes (points) are $x_p=ph$ for $p=1,\dots,m-1$.  This mesh supports a finite-dimensional vector subspace of $\mathcal{H}$:
\begin{equation}
\mathcal{S}^h = \left\{v(x)\,\big|\,v \text{ is continuous, linear on each subinterval, and } v(0)=v(1)=0\right\}.  \label{fespace}
\end{equation}
This space has a basis of ``hat'' functions $\{\psi_p(x)\}$, one for each interior node (Figure \ref{fig:onehat}).  Such a hat function $\psi_p$ is defined by two properties: $\psi_p$ is in $\mathcal{S}^h$ and $\psi_p(x_q)=\delta_{pq}$ for all $q$.  Note that the $L^2$ norm of $\psi_p$ depends on the mesh resolution $h$, and that $\ip{\psi_p}{\psi_q}\ne 0$ for three indices $q=p-1,p,p+1$.  Thus this basis of hat functions, while well-conditioned, is not orthonormal.

The numerical solution $u^h$ has the expansion
\begin{equation}
  u^h(x) = \sum_{p=1}^{m-1} u[p] \psi_p(x)  \label{fesolution}
\end{equation}
with coefficients $u[p]$ equal to the point values $u^h(x_p)$.  That is, because the hat functions form a ``nodal basis'' \cite{Elmanetal2014}, $u^h$ may be represented as a vector $\bu$ in $\RR^{m-1}$ either by its coefficients in the basis $\{\psi_p\}$ or by its point values:
\begin{equation}
\bu =\{u[p]\} = \{u^h(x_p)\}.  \label{fevector}
\end{equation}

The FE approximation $F^h$ of the nonlinear operator $F$ in \eqref{operator} acts on functions in $\mathcal{S}^h$.  Its values $F^h(w^h)[\psi_p]$ are easily computed if the transcendental integral is approximated, for example by using the trapezoid rule, as follows.  Noting that the support of $\psi_p(x)$ is $[x_{p-1},x_{p+1}]$, and that the derivative of $\psi_p$ is $\pm 1/h$, we have:
\begin{align}
  F(w^h)[\psi_p] &= \int_0^1 (w^h)'(x) \psi_p'(x) - \lambda e^{w^h(x)} \psi_p(x)\, dx  \label{feoperator} \\
    &= \int_{x_{p-1}}^{x_{p+1}} (w^h)'(x) (\pm 1/h)\,dx - \lambda \int_{x_{p-1}}^{x_{p+1}} e^{w^h(x)} \psi_p(x)\, dx \notag \\
    &\approx h \left(\frac{w[p]-w[p-1]}{h} - \frac{w[p+1]-w[p]}{h}\right) - h \lambda e^{w[p]}  \notag \\
    &= \frac{1}{h}\left(2w[p]-w[p-1]-w[p+1]\right) - h \lambda e^{w[p]} \notag \\
    &= F^h(w^h)[\psi_p] \notag
\end{align}
Note that $F^h$ is a rescaled version of a well-known $O(h^2)$ finite difference expression.  Function \texttt{FF()} in \texttt{fas1.py} computes this formula.

Now consider the right-hand-side functional $\ell[v]$ in \eqref{weakform}, which we will approximate by $\ell^h[v]$ acting on $\mathcal{S}^h$.  We again apply the trapezoid rule to compute the integral $\ip{g}{\psi_p}$, and we get the simple formula
\begin{equation}
  \ell^h[\psi_p] = h\, g(x_p). \label{ferhs}
\end{equation}
The linear functional $\ell^h$ and the function $g$ are different objects, though they only differ by a factor of the mesh size $h$.

The finite element weak form can now be stated:
\begin{equation}
  F^h(u^h)[v] = \ell^h[v] \qquad \text{for all } v \text{ in } \mathcal{S}^h. \label{feweakform}
\end{equation}
To numerically solve \eqref{feweakform} we will seek an iterate $w^h$ so that the \emph{residual}
\begin{equation}
  r^h(w^h)[v] = \ell^h[v] - F^h(w^h)[v]  \label{feresidual}
\end{equation}
is small for all $v$ in $\mathcal{S}^h$.  Again $r^h(w^h)$ is a linear functional acting on functions in $\mathcal{S}_h$, so it suffices to apply it to a basis of test functions $v=\psi_p$:
\begin{equation}
  r^h(w^h)[\psi_p] = \ell^h[\psi_p] - \frac{1}{h}\left(2w[p]-w[p-1]-w[p+1]\right) + h \lambda e^{w[p]}.  \label{feresidualdetail}
\end{equation}
Solving the finite-dimensional nonlinear system, i.e.~the FE approximation of \eqref{weakform}, is equivalent to finding $w^h$ in $\mathcal{S}^h$ so that $r^h(w^h)[\psi_p]=0$ for $p=1,\dots,m-1$.

A function in \texttt{fas1.py} computes \eqref{feresidualdetail} for any source functional $\ell^h$.  On the original mesh, soon to be called the ``fine mesh'', we will use formula \eqref{ferhs}.  However, the FAS algorithm (sections \ref{sec:fastwolevel} and \ref{sec:cycles}) is a systematic way to introduce a new source functional on each coarser mesh.

The function $u^h(x)$ in $\mathcal{S}^h$, equivalently $\bu$ in $\RR^{m-1}$ given by \eqref{fevector}, exactly solves a finite-dimensional nonlinear system \eqref{feweakform}.  In practice, however, at each stage we only possess an iterate $w^h(x)$, for which the ``algebraic error''
\begin{equation}
  e^h = w^h - u^h  \label{feerror}
\end{equation}
might be small.  On the other hand, $u^h$ is not the continuum solution either; the ``discretization error'' $u^h-u$, where $u$ is the exact solution of the continuum problem \eqref{weakform}, is also nonzero in general.  The theory of an FE method will show that discretization error goes to zero as $h\to 0$, at a particular rate determined by the FE space and the smoothness of the continuum problem \cite{Elmanetal2014}, but such a theory assumes we have exactly-solved the finite-dimensional system, i.e.~that we possess $u^h$ itself.  The full ``numerical error'' is the difference $w^h-u$, and we have
\begin{equation}
\|w^h-u\| \le \|w^h-u^h\|+\|u^h-u\|.
\end{equation}
In other words, the numerical error, which we want to control, is bounded by the algebraic error plus the discretization error.

In the \texttt{-mms} case of \texttt{fas1.py}, where the exact solution $u$ of the continuum problem is known, the numerical error norm $\|w^h-u\|$ is computable.  Normally we cannot access $u$ or $u^h$ directly, and only the residual norm $\|r^h(w^h)\|$ is computable, but the norm of the algebraic error is controlled to within a matrix condition number by the residual norm.

\section{The nonlinear Gauss-Seidel iteration}  \label{sec:ngs}

Next we describe an iteration which will, if carried far enough, solve the finite-dimensional nonlinear system \eqref{feweakform} to desired accuracy.  This is the nonlinear Gauss-Seidel (NGS) iteration \cite{Briggsetal2000}, also called Gauss-Seidel-Newton \cite{BrandtLivne2011}.  It updates the iterate $w^h$ by changing each point value $w^h(x_p)$ to make the residual at that point zero.  That is, NGS solves the problem
\begin{equation}
\phi(c) = r^h(w^h + c \psi_p)[\psi_p] = 0  \label{ngspointproblem}
\end{equation}
for a scalar $c$.  Once $c$ is found we update the point value (coefficient):
\begin{equation}
  w^h \leftarrow w^h + c \psi_p,  \label{ngspointupdate}
\end{equation}
equivalently $w[p] \leftarrow w[p] + c$.

In the linear Gauss-Seidel iteration \cite{Greenbaum1997}, $w[p]$ is updated in a certain nodal ordering, using current values $w[q]$ when evaluating the residual in \eqref{ngspointproblem}.  However, as the residual is made zero at one point it is no longer zero at the previous points.  Gauss-Seidel-type methods are called ``successive'' \cite{GraeserKornhuber2009} or ``multiplicative'' \cite{Bueler2021} corrections.  ``Additive'' corrections, of which the Jacobi iteration \cite{Greenbaum1997} is the best known, are also possible, but they are somewhat less efficient.  Our program only runs in serial, and the parallelizability of the Jacobi iteration cannot be exploited.

Solving the scalar problem $\phi(c)=0$ cannot be done exactly for a transcendental problem like \eqref{liouvillebratu}.  Instead we will use a fixed number of Newton iterations \cite[Chapter 4]{Bueler2021} to generate a (scalar) sequence $\{c_k\}$ converging to $c$.  Starting from $c_0=0$ we compute
\begin{equation}
\phi'(c_k)\, s_k = -\phi(c_k),  \qquad  c_{k+1} = c_k + s_k, \label{ngsnewton}
\end{equation}
for $k=0,1,\dots$.  From \eqref{feresidualdetail} we have
\begin{align*}
 \phi(c) &= \ell^h[\psi_p] - \frac{1}{h} \left(2(w[p]+c) - w[p-1] - w[p+1]\right) + h \lambda e^{w[p]+c}, \\
\phi'(c) &= -\frac{2}{h} + h \lambda e^{w[p]+c}.
\end{align*}
The vast majority of the work of our FAS algorithms will be in evaluating these expressions.

The NGS method sweeps through the mesh, zeroing $\phi(c)$ at successive nodes $x_p$, in increasing $p$ order, as in the following pseudocode which modifies $w^h$ in-place:
\begin{pseudo*}
\pr{ngssweep}(w^h,\ell^h,\id{niters}=2)\text{:} \\+
    $r(w^h)[v] := \ell^h[v] - F^h(w^h)[v]$ \\
    for $p=1,\dots,m-1$ \\+
        $\phi(c) := r^h(w^h + c \psi_p)[\psi_p]$ \\
        $c=0$ \\
        for $k=1,\dots,$\id{niters} \\+
            $c \gets c - \phi(c) / \phi'(c)$ \\-
        $w[p] \gets w[p] + c$
\end{pseudo*}
For FAS algorithms (next section) we also define \textsc{ngssweep-back} with decreasing node order $p=m-1,\dots,1$.  Function \texttt{ngssweep()} in \texttt{fas1.py} computes either order.

For a linear differential equation the Gauss-Seidel iteration is known to converge subject to matrix assumptions which correspond to ellipticity of the original problem \cite[for example]{Greenbaum1997}.  We expect that for weak nonlinearities, e.g.~small $\lambda$ in \eqref{liouvillebratu}, our method will converge as a solution method for \eqref{feweakform}, and we will demonstrate that this occurs in practice (section \ref{sec:convergence}).  However, one observes that, after substantial progress in the first few sweeps during which the residual becomes very smooth, NGS stagnates.  Following Brandt \cite{Brandt1977,BrandtLivne2011}, who asserts that such a stalling scheme must be wrong, we adopt the multigrid approach next.

\section{The FAS equation for two levels}  \label{sec:fastwolevel}

The fundamental goal of any multigrid scheme is to do a minimal amount of work (smoothing) on a given mesh and then switch to a less expensive coarser mesh to do the rest of the work.  By transferring (restricting) a version of the problem to the coarser mesh one can nearly solve for the error.  The coarse-mesh approximation of the error is then added-back (prolonged) to correct the solution on the finer mesh.

Thus the multigrid \emph{full approximation storage} (FAS) scheme \cite{Brandt1977,Briggsetal2000} must include the following elements:
\renewcommand{\labelenumi}{(\roman{enumi})}
\begin{enumerate}
\item a hierarchy of meshes, with restriction and prolongation operators between levels,
\item a ``smoother'' for each level, and
\item a meaningful way to transfer the problem to a coarser mesh.
\end{enumerate}

Regarding (i), we describe only two levels at first, but deep mesh hierarchies are used in section \ref{sec:cycles}.  Here our coarse mesh has spacing $2h$ and $m/2$ elements (subintervals); all quantities on the coarse mesh have superscript ``$2h$''.  The program \texttt{fas1.py} only refines by factors of two, but the ideas generalize for other refinement factors.

For (ii), a small fixed number of NGS sweeps is our smoother.  Each sweep, algorithm \textsc{ngssweep} above, is an $O(m)$ operation with a small constant.  The constant is determined by the number of Newton iterations and the expense of evaluating nonlinearities at each point, e.g.~$\lambda e^u$ in \eqref{liouvillebratu}.  A few NGS sweeps produce smooth fine-mesh residual $r^h(w^h)$ and algebraic error $e^h = w^h - u^h$, but these fields are not necessarily small.  Using more sweeps of NGS would eventually make the error small, and solve problem \eqref{feweakform}, but inefficiently in the sense that many sweeps would be needed, generally giving an $O(m^q)$ method for $q\gg 1$.  However, a coarser mesh will see the coarse-mesh interpolant of the fine-mesh residual as less smooth, so then NGS can quickly eliminate a large fraction of the error.  Descending to yet coarser meshes, in a V-cycle as described in section \ref{sec:cycles}, leads to a coarsest mesh on which the error can be essentially eliminated by applying NGS at only a few interior points.  (In the default settings for \texttt{fas1.py}, the coarsest mesh has two subintervals and one interior point.)

For item (iii), what is the coarse-mesh version of the problem?  To derive this equation, which is Brandt's FAS equation \cite{Brandt1977}, we start from the FE weak form \eqref{feweakform}.  The fine-mesh solution $u^h$ is generally unknown.  For an iterate $w^h$ we subtract $F^h(w^h)[v]$ from both sides of \eqref{feweakform} to get a residual \eqref{feresidual} on the right:
\begin{equation}
  F^h(u^h)[v] - F^h(w^h)[v] = r^h(w^h)[v]. \label{fasproto}
\end{equation}
This is not yet the FAS equation, but three key observations apply to \eqref{fasproto}:
\begin{itemize}
\item Both $w^h$ and $r^h(w^h)$ are known and/or computable.
\item If the smoother has been applied then $e^h=w^h-u^h$ and $r^h(w^h)$ are smooth.
\item If $F^h$ were linear in $w^h$ then we could rewrite \eqref{fasproto} in terms of the error:
    $$\qquad\qquad\qquad\qquad F^h(e^h)[v] = -r^h(w^h)[v], \qquad\qquad (\text{\emph{if $F^h$ is linear}})$$
equivalently $A\be=-\br$ in matrix form.
\end{itemize}

Based on these observations, Brandt proposed a new nonlinear equation for the coarse mesh.  It modifies \eqref{fasproto} by replacing terms using restriction operators on the computable quantities and by re-discretizing the nonlinear operator into $F^{2h}$ acting on $\mathcal{S}^{2h}$.  Because the problem is nonlinear we must store the coarse-mesh solution, namely $u^{2h}$ in $\mathcal{S}^{2h}$, not just the error.  Denoting the restriction operators by $R'$ and $R$, which are addressed in the next section, the following is Brandt's FAS equation:
\begin{equation}
  F^{2h}(u^{2h})[v] - F^{2h}(R w^h)[v] = R' (r^h(w^h))[v], \label{faspreequation}
\end{equation}
for all $v$ in $\mathcal{S}^{2h}$.  We can simplify the appearance by defining a linear functional
\begin{equation}
  \ell^{2h}[v] = R' (r^h(w^h))[v] + F^{2h}(R w^h)[v] \label{fasell}
\end{equation}
so \eqref{faspreequation} becomes
\begin{equation}
  F^{2h}(u^{2h})[v] = \ell^{2h}[v]. \label{fasequation}
\end{equation}

The key idea behind the FAS equation \eqref{fasequation}, which has the same form as the fine-mesh weak form \eqref{feweakform}, is that the smoothness of the error and residual have allowed us to accurately coarsen the problem.  Note that if $w^h=u^h$, that is, if $w^h$ is the exact solution to the fine-mesh problem \eqref{feweakform}, then $r^h(w^h)=0$ so $\ell^{2h}$ simplifies to $F^{2h}(R w^h)[v]$, and the solution of \eqref{fasequation} would be $u^{2h} = R w^h$ by well-posedness.

In stating a two-level FAS method we will suppose \eqref{fasequation} is solved exactly, so $u^{2h}$ and the coarse-mesh error $u^{2h}-Rw^h$ are known.  We then update the fine-mesh iterate by adding a fine-mesh version of the error:
\begin{equation}
  w^h \gets w^h + P(u^{2h} - R w^h) \label{fasupdate}
\end{equation}
Here $P$ is a prolongation operator (next section); it extends a function in $\mathcal{S}^{2h}$ to $\mathcal{S}^h$.  Supposing that the smoother and the restriction/prolongation operators $R',R,P$ are all determined, formulas \eqref{fasequation}, \eqref{fasell}, and \eqref{fasupdate} define the following in-place algorithm in which $F^h$ and $F^{2h}$ denote discretizations of $F$:
\label{fastwolevel}
\begin{pseudo*}
\pr{fas-twolevel}(w^h,\ell^h,\id{down}=1,\id{up}=1)\text{:} \\+
    for $j=1,\dots,$\id{down} \\+
        \pr{ngssweep}(w^h,\ell^h) \\-
    $\ell^{2h}[v] := R' (\ell^h-F^h(w^h))[v] + F^{2h}(R w^h)[v]$ \\
    $w^{2h} = \pr{copy}(R w^h)$ \\
    \pr{coarsesolve}(w^{2h},\ell^{2h}) \\
    $w^h \gets w^h + P(w^{2h} - R w^h)$ \\
    for $j=1,\dots,$\id{up} \\+
        \pr{ngssweep-back}(w^h,\ell^h)
\end{pseudo*}
Note that we allow smoothing before and after the coarse-mesh correction.

While it is common in linear multigrid \cite{Briggsetal2000,Bueler2021,Trottenbergetal2001} to apply a direct solver like LU decomposition as the coarse-mesh solver, our problem is nonlinear and no finite-time direct solver is available.  Instead we do a fixed number of NGS sweeps:
\begin{pseudo*}
\pr{coarsesolve}(w,\ell,\id{coarse}=1)\text{:} \\+
    for $j=1,\dots,$\id{coarse} \\+
        \pr{ngssweep}(w,\ell)
\end{pseudo*}

In order to implement FAS we must define the action of operators $R'$, $R$, and $P$ used in \eqref{fasell} and \eqref{fasupdate},which is done next.  In section \ref{sec:cycles} we will define an FAS V-cycle by replacing \textsc{coarsesolve} with recursive application of the FAS coarsening step.

\section{Restriction and prolongation operators} \label{sec:restrictionprolongation}

To explain the two restriction operators $R'$ and $R$ in \eqref{fasequation}, and the prolongation $P$ in \eqref{fasupdate}, first note that functions $w^h$ in $\mathcal{S}^h$ are distinct objects from linear functionals like the residual $r^h(w^h)$.  Denoting such linear functionals by $(\mathcal{S}^h)'$, the three operators are distinguished by their domain and range spaces:
\begin{align}
  R' &: (\mathcal{S}^h)' \to (\mathcal{S}^{2h})', \label{rpoperators} \\
  R  &: \mathcal{S}^h \to \mathcal{S}^{2h}, \notag \\
  P  &: \mathcal{S}^{2h} \to \mathcal{S}^h. \notag
\end{align}

On the other hand, both functions in $\mathcal{S}^h$ and linear functionals in $(\mathcal{S}^h)'$ are representable by vectors in $\RR^{m-1}$.  One stores a function $w^h$ via coefficients $w[p]$ with respect to an expansion in the hat function basis $\{\psi_p\}$, as in \eqref{fesolution} for example, while one stores a functional $\ell^h$ by its values $\ell^h[\psi_p]$.  Though it makes sense to represent $w^h$ as a column vector and $\ell^h$ as a row vector \cite{TrefethenBau1997}, in Python we will use ``flat'' one-dimensional NumPy arrays \cite{Harrisetal2020} for both purposes.  For our problem an iterate $w^h$ has zero boundary values, and likewise $\ell^h$ acts on functions having zero boundary values, thus only interior hat functions are needed in these representations.

\begin{figure}
\includegraphics[width=0.6\textwidth]{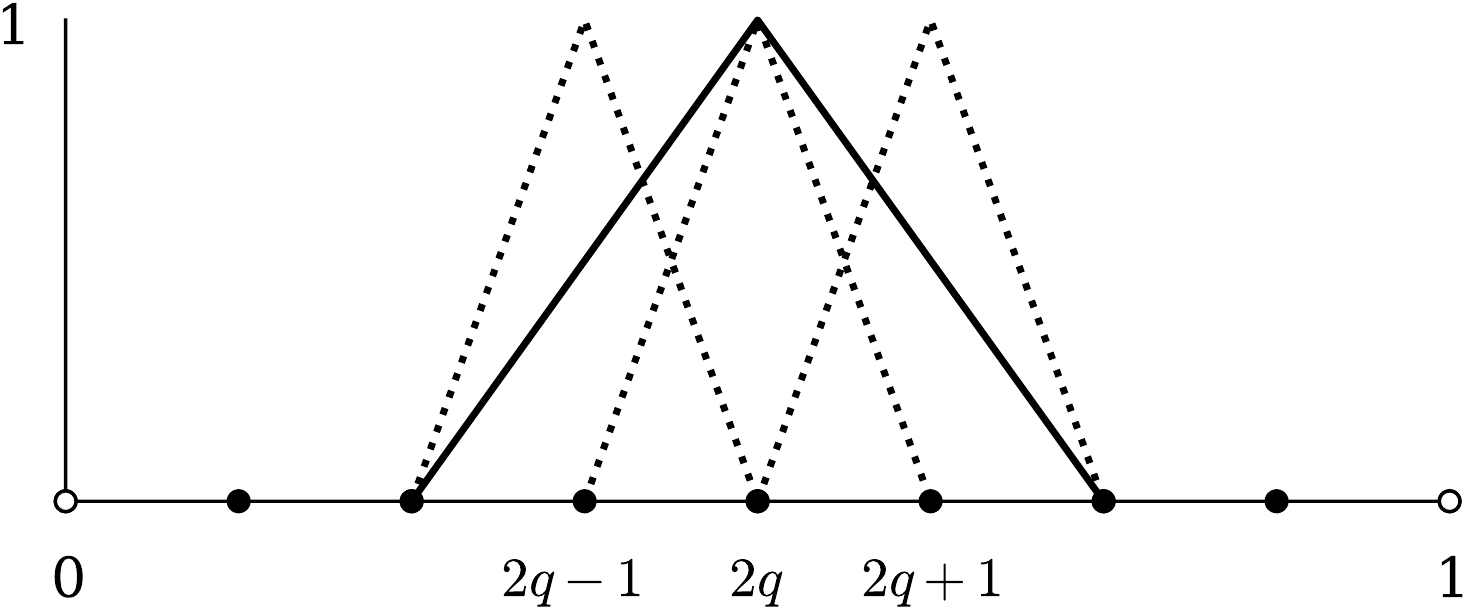}
\caption{Formula \eqref{hatrelation} writes a coarse-mesh hat function $\psi_q^{2h}(x)$ (solid) as a linear combination of fine-mesh hats $\psi_p^h(x)$ (dotted).}
\label{fig:hatcombination}
\end{figure}

But how do $R'$, $R$, and $P$ actually operate?  The key calculation relates the coarse-mesh hat functions $\psi_q^{2h}(x)$ to the fine mesh hats $\psi_p^h(x)$:
\begin{equation}
  \psi_q^{2h}(x) = \frac{1}{2} \psi_{2q-1}^h(x) + \psi_{2q}^h(x) + \frac{1}{2} \psi_{2q+1}^h(x), \label{hatrelation}
\end{equation}
for $q=1,2,\dots,m/2-1$.  (We are assuming $m$ is even.)  See Figure \ref{fig:hatcombination}.

First consider the prolongation $P$.  Because a piecewise-linear function on the coarse mesh is also a piecewise-linear function on the fine mesh, in this FE context $P$ is defined as the injection of $\mathcal{S}^{2h}$ into $\mathcal{S}^h$, without changing the function.  If $w^{2h}(x) = \sum_{q=1}^{m/2-1} w[q] \psi_q^{2h}(x)$ then \eqref{hatrelation} computes $P w^{2h}$ in terms of fine-mesh hat functions:
\begin{align}
(P w^{2h})(x) &= \sum_{q=1}^{m/2-1} w[q] \left(\frac{1}{2} \psi_{2q-1}^h(x) + \psi_{2q}^h(x) + \frac{1}{2} \psi_{2q+1}^h(x)\right) \label{pformula} \\
              &= \frac{1}{2} w[1] \psi_1^h(x) + w[1] \psi_2^h(x) + \left(\frac{1}{2} w[1] + \frac{1}{2} w[2]\right) \psi_3^h(x) + w[2] \psi_4^h(x) \notag \\
              &\qquad + \left(\frac{1}{2} w[2] + \frac{1}{2} w[3]\right) \psi_5^h(x) + \dots + w[m/2\!-\!1] \psi_{m-2}^h(x) \notag \\
              &\qquad + \frac{1}{2} w[m/2\!-\!1] \psi_{m-1}^h(x) \notag
\end{align}
As a matrix, $P:\RR^{m/2-1} \to \RR^{m-1}$ acts on vectors; it has $m/2-1$ columns and $m-1$ rows:
\begin{equation}
P = \begin{bmatrix}
1/2 & & & \\
1 & & & \\
1/2 & 1/2 & & \\
 & 1 & & \\
 & 1/2 & 1/2 & \\
 & & & \ddots
\end{bmatrix} \label{pmatrix}
\end{equation}
The columns of $P$ are linearly-independent and sum to 2 by \eqref{hatrelation}.  The row sums equal one except for the first and last rows.

Next, the restriction $R'$ acts on fine-mesh linear functionals $\ell:\mathcal{S}^h \to \RR$.  It is called ``canonical restriction'' \cite{GraeserKornhuber2009} because the output functional $R'\ell:\mathcal{S}^{2h}\to \RR$ acts on coarse-mesh functions in the same way $\ell$ itself acts on those functions.  Defining $R'$ involves no choices.  We may state this using $P$: for $v$ in $\mathcal{S}^{2h}$,
\begin{equation}
  (R'\ell)[v] = \ell[Pv].  \label{rprimedefinition}
\end{equation}
As noted earlier, $\ell$ in $\RR^{m-1}$ has entries $\ell[\psi_p^h]$.  One computes the values of $R'\ell$ using \eqref{hatrelation}:
\begin{align}
  (R'\ell)[\psi_q^{2h}] &= \ell[\psi_q^{2h}] = \ell\left[\frac{1}{2} \psi_{2q-1}^h + \psi_{2q}^h + \frac{1}{2} \psi_{2q+1}^h\right]  \label{rprimeformula} \\
      &= \frac{1}{2} \ell[\psi_{2q-1}^h] + \ell[\psi_{2q}^h] + \frac{1}{2} \ell[\psi_{2q+1}^h].  \notag
\end{align}
As a matrix, $R'$ is the transpose of $P$, with $m/2-1$ rows and $m-1$ columns:
\begin{equation}
R' = \begin{bmatrix}
1/2 & 1 & 1/2 &   &     & \\
    &   & 1/2 & 1 & 1/2 & \\
    &   &     &   & 1/2 & \\
    &   &     &   &     & \ddots
\end{bmatrix} \label{rprimematrix}
\end{equation}

Finally we consider the restriction $R:\mathcal{S}^h\to\mathcal{S}^{2h}$ acting on functions, a more interesting map which loses information.  (By contrast, $P$ and $R'$ preserve their input, without loss, via reinterpretation on the output mesh.)  Consider a fine-mesh function $w^h = \sum_{p=1}^{m-1} w[p] \psi_p^{h}$.  The result $R w^h$ should be linear across those fine-mesh nodes which are not in the coarse mesh, so values at in-between nodes are not recoverable.

There are three well-known choices for the restriction $R$:
\begin{itemize}
\item $\Rpr$ is defined as projection, by the property
\begin{equation}
  \ip{\Rpr w^h}{v} = \ip{w^h}{v} \label{rprdefinition}
\end{equation}
for all $v\in \mathcal{S}^{2h}$.  Computing the entries of $\Rpr$ requires solving a linear system.  To do so we define invertible, sparse, symmetric mass matrices $Q_{jk}^{h} = \ip{\psi_j^{h}}{\psi_k^{h}}$ for the fine mesh and $Q_{jk}^{2h} = \ip{\psi_j^{2h}}{\psi_k^{2h}}$ for the coarse \cite{Elmanetal2014}.  Then one solves a matrix equation for $\Rpr$:
\begin{equation}
  Q^{2h} \Rpr = R' Q^{h},  \label{rprequation}
\end{equation}
or equivalently $\Rpr = (Q^{2h})^{-1} R' Q^{h}$.  Equation \eqref{rprequation} is justified by using $v=\psi_s^{2h}$ in definition \eqref{rprdefinition}, and then applying \eqref{hatrelation}, as follows.  Write $z = \Rpr w^h = \sum_{q=1}^{M-1} z[q] \psi_q^{2h}$ and expand both sides:
\begin{align*}
\ip{z}{\psi_s^{2h}} &= \ip{w^h}{\psi_s^{2h}} \\
\sum_{q=1}^{m/2-1} z[q] \ip{\psi_q^{2h}}{\psi_s^{2h}} &= \sum_{p=1}^{m-1} w[p] \ip{\psi_p^{h}}{\frac{1}{2} \psi_{2s-1}^{h} + \psi_{2s}^{h} + \frac{1}{2} \psi_{2s+1}^{h}} \\
\sum_{q=1}^{m/2-1} Q_{sq}^{2h} z[q] &= \sum_{p=1}^{m-1} \left(\frac{1}{2} Q_{2s-1,p} + Q_{2s,p} + \frac{1}{2} Q_{2s+1,p}\right) w[p] \\
(Q^{2h} \Rpr w^h)[s] &= (R' Q^h w^h)[s]
\end{align*}
(Note $w^h$ in $\mathcal{S}^h$ and index $s$ are arbitrary.)  In 1D the mass matrices $Q^{2h},Q^h$ are tridiagonal, thus each column of $\Rpr$ can be found by solving equation \eqref{rprequation} using an $O(m)$ algorithm \cite{TrefethenBau1997}.  While this is implementable, and computable by hand in this case, the alternatives below are easier to implement.
\item $\Rin$ is defined as pointwise injection.  Supposing $w^h = \sum_{p=1}^{m-1} w[p] \psi_p^{h}$,
\begin{equation}
  \Rin w^h = \sum_{q=1}^{m/2-1} w[2q] \psi_q^{2h}, \label{rindefinition}
\end{equation}
so $(\Rin w^h)(x_q) = w^h(x_q) = w[2q]$ for each point $x_q$.  In other words, to compute $\Rin w^h$ we simply drop the nodal values at those fine-mesh nodes which are not in the coarse mesh.  As a matrix this is
\begin{equation}
\Rin = \begin{bmatrix}
0 & 1 &   &   &   &   &\\
  &   & 0 & 1 &   &   & \\
  &   &   &   & 0 & 1 & \\
  &   &   &   &   &   & \ddots
\end{bmatrix}. \label{rinmatrix}
\end{equation}
This restriction is very simple but it can lose track of the magnitude of $w^h$, or badly mis-represent it, \emph{if} the input is not smooth.  For example, sampling a sawtooth function at the coarse-mesh nodes would capture only the peaks or only the troughs.
\item $\Rfw$, the ``full-weighting'' restriction \cite{Briggsetal2000}, averages nodal values onto the coarse mesh:
\begin{equation}
  \Rfw w^h = \sum_{q=1}^{m/2-1} \left(\frac{1}{4} w[2q-1] + \frac{1}{2} w[2q] + \frac{1}{4} w[2q+1]\right) \psi_q^{2h}. \label{rfwdefinition}
\end{equation}
This computes each coarse-mesh nodal value of $z=\Rfw w^h$ as a weighted average of the value of $w^h$ at the three closest fine-mesh nodes.  The matrix is a multiple of the canonical restriction matrix in \eqref{rprimematrix}:
\begin{equation}
\Rfw = \begin{bmatrix}
1/4 & 1/2 & 1/4 &     &     &  \\
    &     & 1/4 & 1/2 & 1/4 &  \\
    &     &     &     & 1/4 &  \\
    &     &     &     &     & \ddots
\end{bmatrix} = \frac{1}{2} R'. \label{rfwmatrix}
\end{equation}
\end{itemize}

\medskip
Which restriction do we choose?  Because of their simplicity, we will implement and compare $\Rin$ and $\Rfw$ in \texttt{fas1.py}.

\section{Cycles} \label{sec:cycles}

The main principles of the FAS scheme are already contained in the \textsc{fas-twolevel} algorithm in section \ref{sec:fastwolevel}, from which it is a small step to solve the coarse-mesh problem by the same scheme, creating a so-called ``V-cycle''.  To define this precisely we need an indexed hierarchy of mesh levels.  Start with a coarsest mesh with $m_0$ elements of length $h_0=1/m_0$.  (By default in \texttt{fas1.py} we have $m_0=2$.)  For $k=1,\dots,K$ we refine by factors of two so that the $k$th mesh has $m_k=2^k m_0$ elements of length $h_k=h_0/2^k$.  The final $K$th mesh is now called the ``fine mesh''.  Instead of the superscripts $h$ and $2h$ used in section \ref{sec:fastwolevel}, a ``$k$'' superscript indicates the mesh on which a quantity lives.

On this hierarchy an FAS V-cycle is the following in-place recursive algorithm:
\begin{pseudo*}
\pr{fas-vcycle}(k,w^k,\ell^k,\id{down}=1,\id{up}=1)\text{:} \\+
    if $k=0$ \\+
        \pr{coarsesolve}(w^0,\ell^0) \\-
    else \\+
        for $j=1,\dots,$\id{down} \\+
            \pr{ngssweep}(w^k,\ell^k) \\-
        $w^{k-1} = \pr{copy}(R w^k)$ \\
        $\ell^{k-1}[v] := R' (\ell^k-F^k(w^k))[v] + F^{k-1}(R w^k)[v]$ \\
        \pr{fas-vcycle}(k-1,w^{k-1},\ell^{k-1}) \\
        $w^k \gets w^k + P(w^{k-1} - R w^k)$ \\
        for $j=1,\dots,$\id{up} \\+
            \pr{ngssweep-back}(w^k,\ell^k) \\-
\end{pseudo*}

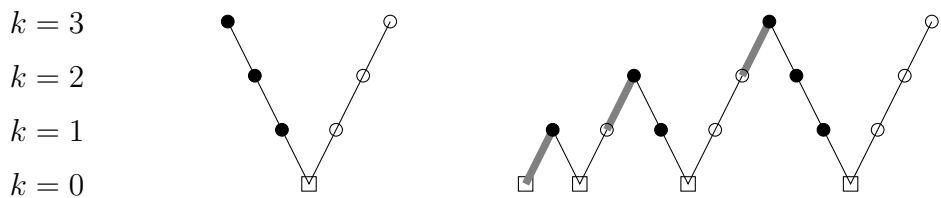
\begin{figure}
\begin{tikzpicture}[scale=1.2]
  \pgfmathsetmacro\hstep{0.3}
  \pgfmathsetmacro\vstep{0.6}
  \pgfmathsetmacro\ceps{0.08}   

  \node at (-2,3*\vstep) {$k=3$};
  \node at (-2,2*\vstep) {$k=2$};
  \node at (-2,\vstep) {$k=1$};
  \node at (-2,0.0) {$k=0$};

  \draw[black,thin] (0.0,3*\vstep) -- (\hstep,2*\vstep) --  (2*\hstep,\vstep) -- (3*\hstep,0.0)
                    -- (4*\hstep,\vstep) -- (5*\hstep,2*\vstep) -- (6*\hstep,3*\vstep);
  \filldraw (0.0,3*\vstep) circle (2.0pt);
  \filldraw (\hstep,2*\vstep) circle (2.0pt);
  \filldraw (2*\hstep,\vstep) circle (2.0pt);
  \draw     (3*\hstep-\ceps,-\ceps) rectangle (3*\hstep+\ceps,+\ceps);
  \draw     (4*\hstep,\vstep) circle (2.0pt);
  \draw     (5*\hstep,2*\vstep) circle (2.0pt);
  \draw     (6*\hstep,3*\vstep) circle (2.0pt);

  \pgfmathsetmacro\hoff{11*\hstep}
  \draw[shift={(\hoff,0)}]     (-\ceps,-\ceps) rectangle (+\ceps,+\ceps);
  \draw[line width=1mm,gray,shift={(\hoff,0)}] (0.0,0.0) -- (\hstep,\vstep);

  \pgfmathsetmacro\hoff{12*\hstep}
  \draw[shift={(\hoff,0)},black,thin] (0.0,\vstep) -- (\hstep,0.0) -- (2*\hstep,\vstep);
  \draw[line width=1mm,gray,shift={(\hoff,0)}] (2*\hstep,\vstep) -- (3*\hstep,2*\vstep);
  \filldraw[shift={(\hoff,0)}] (0.0,\vstep) circle (2.0pt);
  \draw[shift={(\hoff,0)}]     (\hstep-\ceps,-\ceps) rectangle (\hstep+\ceps,+\ceps);
  \draw[shift={(\hoff,0)}]     (2*\hstep,\vstep) circle (2.0pt);

  \pgfmathsetmacro\hoff{15*\hstep}
  \draw[shift={(\hoff,0)},black,thin] (0.0,2*\vstep) --  (\hstep,\vstep) -- (2*\hstep,0.0) -- (3*\hstep,\vstep) -- (4*\hstep,2*\vstep);
  \draw[line width=1mm,gray,shift={(\hoff,0)}] (4*\hstep,2*\vstep) -- (5*\hstep,3*\vstep);
  \filldraw[shift={(\hoff,0)}] (0.0,2*\vstep) circle (2.0pt);
  \filldraw[shift={(\hoff,0)}] (\hstep,\vstep) circle (2.0pt);
  \draw[shift={(\hoff,0)}]     (2*\hstep-\ceps,-\ceps) rectangle (2*\hstep+\ceps,+\ceps);
  \draw[shift={(\hoff,0)}]     (3*\hstep,\vstep) circle (2.0pt);
  \draw[shift={(\hoff,0)}]     (4*\hstep,2*\vstep) circle (2.0pt);

  \pgfmathsetmacro\hoff{20*\hstep}
  \draw[shift={(\hoff,0)},black,thin] (0.0,3*\vstep) -- (\hstep,2*\vstep) --  (2*\hstep,\vstep) -- (3*\hstep,0.0)
                    -- (4*\hstep,\vstep) -- (5*\hstep,2*\vstep) -- (6*\hstep,3*\vstep);
  \filldraw[shift={(\hoff,0)}] (0.0,3*\vstep) circle (2.0pt);
  \filldraw[shift={(\hoff,0)}] (\hstep,2*\vstep) circle (2.0pt);
  \filldraw[shift={(\hoff,0)}] (2*\hstep,\vstep) circle (2.0pt);
  \draw[shift={(\hoff,0)}]     (3*\hstep-\ceps,-\ceps) rectangle (3*\hstep+\ceps,+\ceps);
  \draw[shift={(\hoff,0)}]     (4*\hstep,\vstep) circle (2.0pt);
  \draw[shift={(\hoff,0)}]     (5*\hstep,2*\vstep) circle (2.0pt);
  \draw[shift={(\hoff,0)}]     (6*\hstep,3*\vstep) circle (2.0pt);
\end{tikzpicture}
\caption{An FAS V-cycle (left) and F-cycle (right) on a mesh hierarchy with four levels ($K=3$).  Solid dots are \texttt{down} sweeps of NGS, open circles are \texttt{up} sweeps, and squares are \textsc{coarsesolve}. Thick grey edges show $\hat P$.}
\label{fig:cycles}
\end{figure}

The definition of $\ell^k$ depends on the mesh level.  On the fine level it is $\ell^K[v] = \ip{g}{v}$, as in \eqref{ferhs}, but on coarser levels it is determined by the nontrivial FAS formula \eqref{fasell}.  A V-cycle with $K=3$ is shown in Figure \ref{fig:cycles}.  V-cycles can be iterated to solve problem \eqref{feweakform} to desired accuracy:
\begin{pseudo*}
\pr{fas-solver}(w^K,\id{rtol}=10^{-4},\id{cyclemax}=100)\text{:} \\+
    $\ell^K[v] = \ip{g}{v}$ \\
    $r_0 = \|\ell^K - F^K(w^K)\|$ \\
    for $s=1,\dots,\id{cyclemax}$ \\+
        \pr{fas-vcycle}(K,w^K,\ell^K) \\
        if $\|\ell^K-F^K(w^K)\| < \id{rtol}\,r_0$ \\+
            break \\--
    return $w^K$
\end{pseudo*}

Our Python code \texttt{fas1.py} implements \pr{fas-vcycle} and \pr{fas-solver}.  Options \texttt{-rtol}, \texttt{-cyclemax} override the defaults for \pr{fas-solver}.  As demonstrated in the next two sections, 7 to 12 V-cycles, using the default settings in \textsc{fas-vcycle}, with \texttt{down} $=1$ and \texttt{up} $=1$ smoother applications, make a very effective solver on any mesh.

However, we can add a different multilevel idea to get an even better cycle.  It is based on the observation that an iterative equation solver, linear or nonlinear, often depends critically on the quality of its initial iterate.  Indeed, choosing initial iterate $w^K=0$ and calling \textsc{fas-solver} may not yield a convergent method.  However, one finds in practice that coarse meshes are more forgiving with respect to the initial iterate than are finer meshes.  The new idea is to start on the coarsest mesh in the hierarchy, where a blind guess like $w^0=0$ is most likely to succeed, and then work upward through the levels.  At each mesh level one computes an initial iterate by prolongation of a nearly-converged iterate on the previous level, and then one does a V-cycle.  At the finest mesh level we may do repeated V-cycles.

The resulting algorithm is called an FAS multigrid ``F-cycle'' because the pattern in Figure \ref{fig:cycles} (right) looks vaguely like an ``F'' on its back:
\begin{pseudo*}
\pr{fas-fcycle}(K)\text{:} \\+
    $w^0 = 0$ \\
    $\ell^0[v] = \ip{g}{v}$ \\
    \pr{coarsesolve}(w^0,\ell^0) \\
    for $k=1,\dots,K$ \\+
        $w^k = \hat P w^{k-1}$ \\
        $\ell^k[v] = \ip{g}{v}$ \\
        \pr{fas-vcycle}(k,w^k,\ell^k) \\-
    return $w^K$
\end{pseudo*}

This algorithm is also called a ``full multigrid'' (FMG) cycle \cite{BrandtLivne2011,Briggsetal2000}, but the meaning of ``full'' is fundamentally different in FAS versus FMG terminology.

One may run \pr{fas-fcycle} to generate the initial iterate for \pr{fas-solver}.  However, as seen in section \ref{sec:performance} the result of one F-cycle is already a very good solution.

It is important to avoid the introduction of high frequencies as one generates the first iterate on the finer mesh.  Thus a coarse-mesh solution is prolonged on to the next level by a possibly-different operator:
\begin{equation}
  w^k = \hat P w^{k-1} \label{enhancedprolongation}
\end{equation}
It is common for a better interpolation scheme to be used for $\hat P$ than for $P$ \cite{Trottenbergetal2001}.  Our choice for $\hat P$ first applies $P$ to generate a fine-mesh function, but this is followed by sweeping once through the \emph{new} fine-mesh nodes and applying NGS there without altering values at the coarse-mesh nodes.  This $\hat P$ is half of a smoother, and counted as such; see section \ref{sec:performance}.

\section{Demonstrations, and convergence} \label{sec:convergence}

The Python program \texttt{fas1.py} accompanying this note applies \pr{fas-solver} by default, with zero initial iterate, to solve equation \eqref{liouvillebratu}.  The program depends only on the widely-available NumPy library \cite{Harrisetal2020}.

To get started, clone the Git repository and run the program:
\begin{cline}
$ git clone https://github.com/bueler/fas-intro.git
$ cd fas-intro/fas/py/
$ ./fas1.py
  m=8 mesh, 6 V(1,1) cycles (19.50 WU): |u|_2=0.102443
\end{cline}
The V-cycles in this run (Figure \ref{fig:cycles}) are reported as ``\texttt{V(1,1)}'' because the defaults correspond to \texttt{down} $=1$ and \texttt{up} $=1$ NGS sweeps on each level.

Various allowed options are shown by usage help:
\begin{cline}
$ ./fas1.py -h
\end{cline}
Also, a small suite of software (regression) tests of \texttt{fas1.py} is run via \,\texttt{make test}.

\begin{figure}
\includegraphics[width=0.8\textwidth]{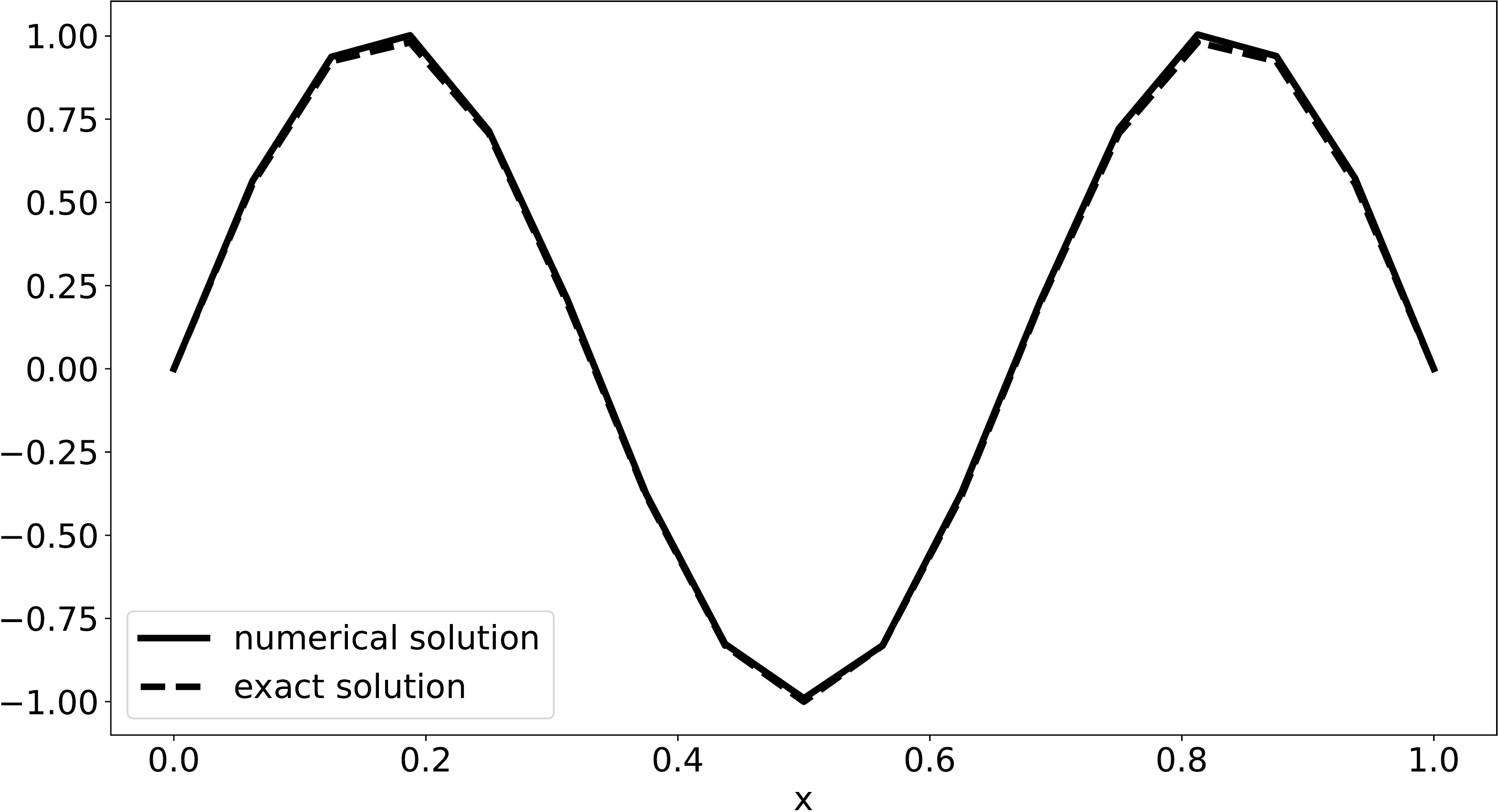}
\caption{Results from a \texttt{-mms} run of \texttt{fas1.py} on $m=16$ elements.}
\label{fig:show}
\end{figure}

Choosing a mesh with $m=2^{K+1}=16$ elements and a problem with known exact solution (section \ref{sec:intro}) yields Figure \ref{fig:show}:
\begin{cline}
$ ./fas1.py -K 3 -mms -show
  m=16 mesh, 6 V(1,1) cycles (21.75 WU): ... |u-u_ex|_2=2.1315e-02
\end{cline}
Note that runs with option \texttt{-mms} report the final numerical error $\|w^h-u\|_2$.

\begin{figure}
\includegraphics[width=0.65\textwidth]{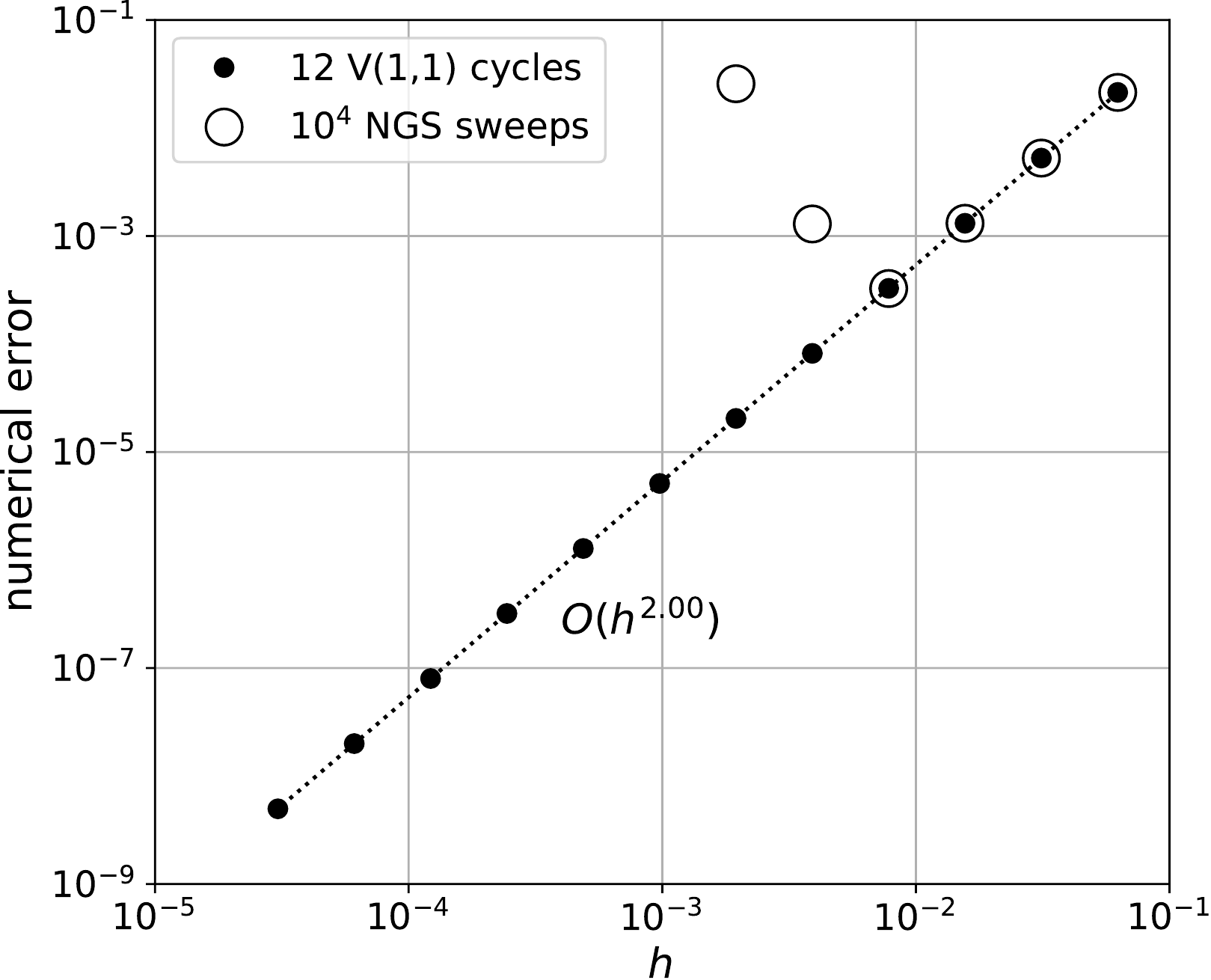}
\caption{For a fixed number of V-cycles the numerical error $\|u-u_{\text{ex}}\|_2$ converges to zero at the expected rate $O(h^2)$.  Even $10^4$ NGS sweeps fail to converge at higher resolutions.}
\label{fig:converge}
\end{figure}

We can demonstrate convergence of our implemented FE method, and verify \texttt{fas1.py}.  Consider hierarchy depths $K=3,4,\dots,14$, corresponding to meshes with $16\le m \le 32768$ elements.  The numerical errors from runs with 12 V-cycles, i.e.~options \texttt{-rtol 0 -cyclemax 12}, are shown in Figure \ref{fig:converge}.  Because our problem is so simple, with a very smooth solution, the convergence rate is exactly at the expected rate $O(h^2)$ \cite{Elmanetal2014}.

However, if instead of a small, fixed number of V-cycles we instead try a large number of NGS sweeps, e.g.~we apply the algorithm below with \texttt{-rtol 0 -cyclemax 10000} and zero initial iterate, then the performance is far to slow to repeat these verification runs.
\begin{pseudo*}
\pr{ngsonly}(w^K,\id{rtol}=10^{-4},\id{cyclemax}=100)\text{:} \\+
    $\ell^K[v] = \ip{g}{v}$ \\
    $r_0 = \|\ell^K - F^K(w^K)\|$ \\
    for $s=1,\dots,\id{cyclemax}$ \\+
        \pr{ngssweep}(w^K,\ell^K) \\
        if $\|\ell^K-F^K(w^K)\| < \id{rtol}\,r_0$ \\+
            break \\--
    return $w^K$
\end{pseudo*}
As shown in Figure \ref{fig:converge}, such runs generates convergence to discretization error only on the 4 coarsest meshes.  For slightly finer meshes ($m=256,512$) the same number of sweeps is no longer sufficient, and on yet finer meshes the same number of sweeps make essentially no progress (not shown).

The reason for the failure of \pr{ngsonly} is that almost all of the algebraic error (section \ref{sec:femethod}) is in low-frequency modes which the NGS sweeps are barely able to reduce.  This is exactly the situation which multigrid schemes are designed to address \cite{BrandtLivne2011,Briggsetal2000}: by moving the problem between meshes the same smoother will efficiently-reduce all frequencies present in the error.  Both the smoother and the coarse-level solver components of our FAS algorithms consist entirely of NGS sweeps, but by adding a multilevel mesh infrastructure we have arranged that the sweeps are always making progress.

\section{Performance}  \label{sec:performance}

Having verified our method, our first performance test compares three solver algorithms:
\begin{itemize}
\item \textsc{fas-fcycle}, defined in section \ref{sec:cycles}.
\item \textsc{fas-solver}, which does V-cycles, also defined in section \ref{sec:cycles}.
\item \textsc{ngsonly} in section \ref{sec:convergence}.
\end{itemize}
The two FAS pseudocodes actually represent many different algorithms according to the different options.  While making no attempt to systematically-explore the parameter space, we observe that 7 to 12 V(1,1) cycles suffice to approach discretization error in the \texttt{-mms} problem.  For F-cycles we must choose how many V-cycles to take once the finest level is reached, and 2 or 3 usually suffice.  Experimentation in minimizing the work units (below), while maintaining convergence, yields a choice of three V(1,0) cycles.

Consider the following specific \texttt{fas1.py} options on meshes with $m=2^{K+1}$ elements for $K=3,4,\dots,17,18$:

\medskip
\begin{tabular}{ll}
\textsf{F-cycle$+$3$\times$V(1,0)} \,:        &\texttt{-mms -fcycle -rtol 0 -cyclemax 4 -up 0 -K }$K$ \\
\textsf{12 V(1,1) cycles} \,:  &\texttt{-mms -rtol 0 -cyclemax 12 -K }$K$ \\
\textsf{NGS sweeps} \,:      &\texttt{-mms -rtol 0 -cyclemax $Z$ -ngsonly -K }$K$
\end{tabular}

\medskip
In order to achieve convergence for NGS sweeps alone, we must choose rapidly increasing $Z$ as $K$ increases.  For the comparison below we simply double $Z$ until the reported numerical error is within a factor of two of discretization error (as reported by the FAS algorithms).  However, at $K=7$ the time is 100 seconds and we stop testing.

\begin{figure}
\includegraphics[width=0.65\textwidth]{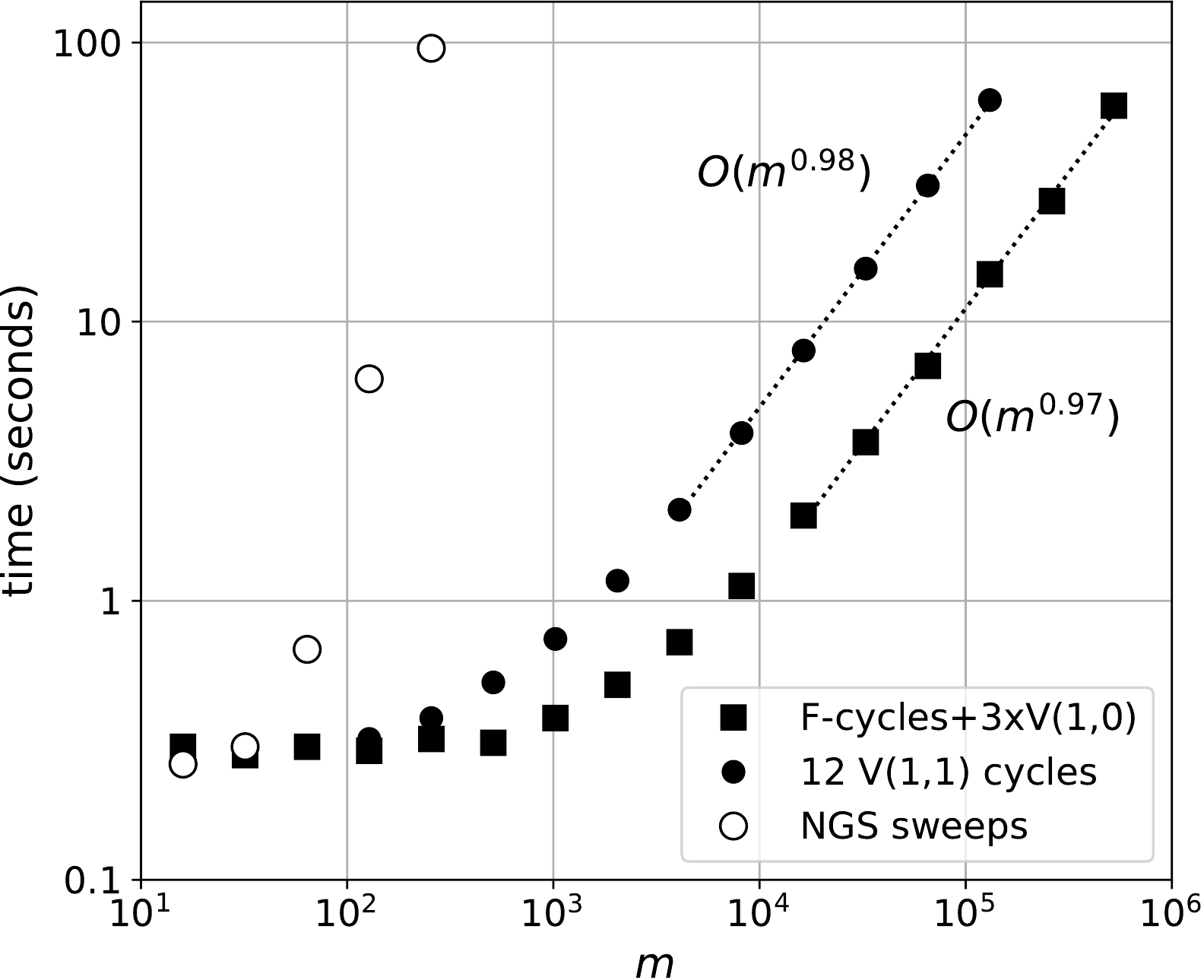}
\caption{Run time to reach discretization error is optimal $O(m)$ for both V-cycles and F-cycles.  Run time explodes for NGS sweeps.}
\label{fig:optimal}
\end{figure}

Run times on the author's laptop are shown in Figure \ref{fig:optimal}.  For all the coarser meshes, e.g.~$m=16,\dots,256$, the FAS algorithms run in about 0.3 seconds.  This is the minimum time to start and run any Python program on this machine, so the actual computational time is not actually observed.  For $m \ge 10^3$ both FAS algorithms enter into a regime where the run time is greater than one second, and then it becomes proportional to $m$.  That is, their solver complexity is $O(m^1)$.  These are \emph{optimal} solvers \cite[Chapter 7]{Bueler2021}.  By contrast, \pr{ngsonly} is far from optimal, and indeed it is not capable of solving on fine meshes.  Fitting the three finest-mesh completed cases suggests its time is $O(m^{3.5})$.

A standard way to compare multigrid-type solvers uses the concept of a \emph{work unit} (WU).  One WU is the number of operations needed to do a smoother sweep on the finest mesh, which takes $O(m)$ arithmetic (floating point) operations.  For WUs in a 1D multilevel scheme, a smoother sweep on the second-finest mesh is $\frac{1}{2}$WU, and so on downward in the hierarchy.  The total of WU for a multigrid algorithm is thus a finite geometric sum \cite{Briggsetal2000} which depends on the number of levels $K$.  For simplicity we do not count the arithmetic work in restriction and prolongation, other than in the enhanced prolongation $\hat P$ in \eqref{enhancedprolongation}, which uses $\frac{1}{2}$WU when passing to the finest mesh.  Also we ignore non-arithmetic work entirely, for example vector copies.

Consider the $K\to\infty$ limit of WU calculations for the three algorithms above:
\begin{align*}
\text{WU}\big(\text{\textsf{F-cycle$+Z\times$V(1,0)}}\big) &\approx 3+2Z \\
\text{WU}\big(\text{\textsf{$Z$ V(1,1) cycles}}\big)   &\approx 4Z \\
\text{WU}\big(\text{\textsf{$Z$ NGS sweeps}}\big)      &= Z
\end{align*}
To confirm these estimates in practice we have added WU counting in \texttt{fas1.py}.  On a $K=10$ mesh with $m=2^{11}=2048$ elements, for example, we observe that \textsf{F-cycles$+$3$\times$V(1,0)} requires a measured 8.98 WU while \textsf{12 V(1,1) cycles} uses 47.96 WU.

A single F-cycle, without any additional V-cycles, nearly reaches discretization error.  Consider three single-F-cycle schemes.  The first is ``F(1,1)'', which uses the default settings \id{down}=1 and \id{up}=1.  The other two are ``F(1,0)'', using \id{up}=0, and ``F(1,0)+$\Rin$'', which changes from the default full-weighting restriction ($\Rfw$) to injection ($\Rin$).  These three solvers use 9, 5, and 5 WU, respectively, in the $K\to\infty$ limit of many levels.

\begin{figure}
\includegraphics[width=0.65\textwidth]{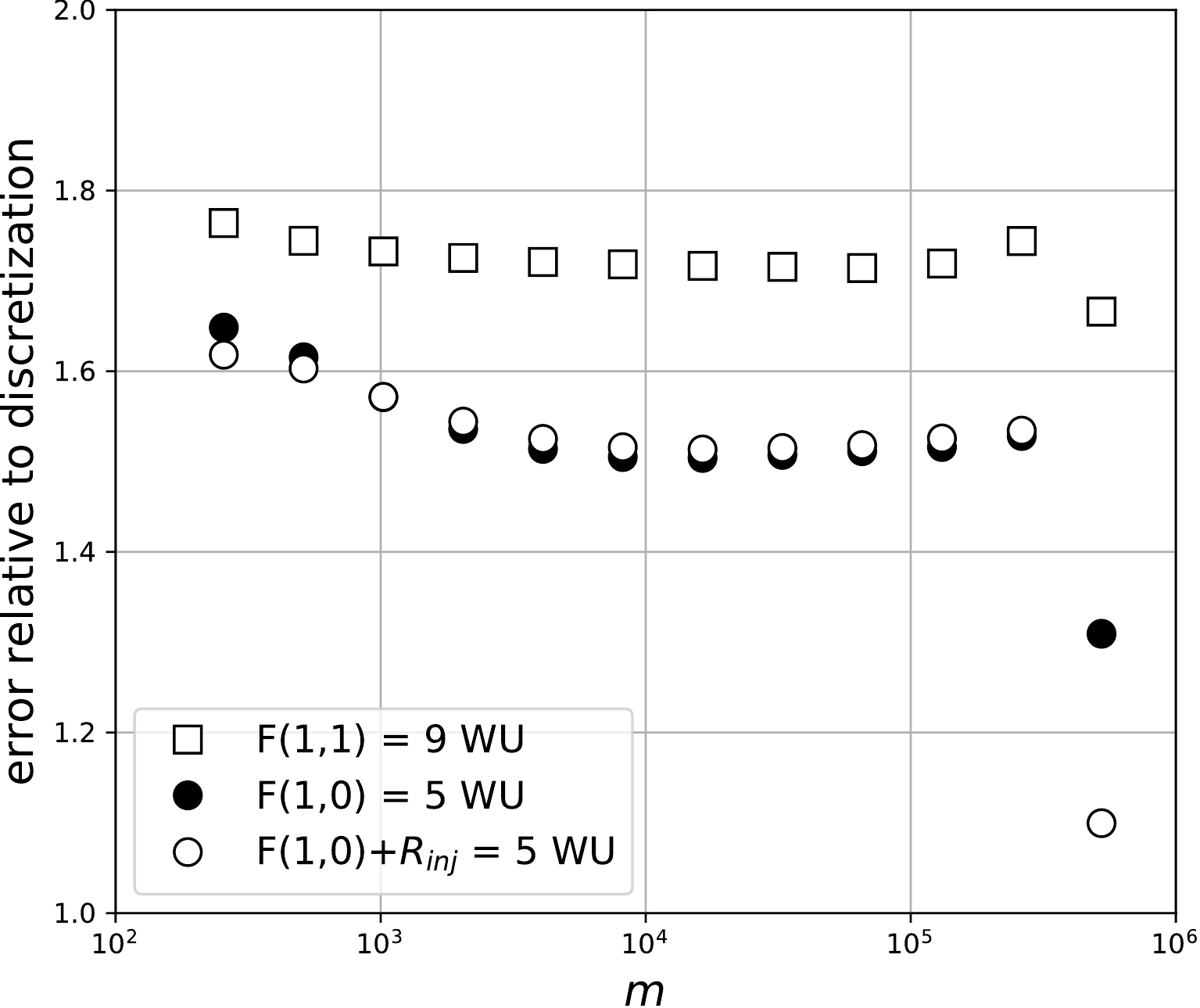}
\caption{Computed numerical error, relative to discretization error, from three versions of a single F-cycle.}
\label{fig:tme}
\end{figure}

Figure \ref{fig:tme} shows that on $K=7,\dots,18$ meshes, with up to $m=2^{19} = 5 \times 10^5$ elements, the measured numerical error is within a factor of two of discretization error.  Note that the F(1,0) cycles actually generate smaller errors.  There is no significant difference between the two restriction methods.  On the finest mesh it seems that the discretization error itself, of order $10^{-11}$, was corrupted by rounding errors.  Noting that some multigrid authors \cite[for example]{BrownSmithAhmadia2013} assert ``textbook multigrid efficiency'' for a scheme when fewer than 10 WU are needed to achieve discretization error, we conclude that our F-cycles exhibit textbook multigrid efficiency.

\section{Extensions}  \label{sec:extensions}

Our program \texttt{fas1.py} is deliberately basic in many senses.  Here are three possible extensions which the reader might want to implement:
\renewcommand{\labelenumi}{\textbf{\Roman{enumi}.}}
\begin{enumerate}
\item The default value of the parameter $\lambda$ in \eqref{liouvillebratu} is \texttt{-lam 1.0}, but one can check that the $g=0$ problem becomes unstable at a critical value $\lambda_c \approx 3.5$.  Interestingly, the solution changes very little as $\lambda \nearrow \lambda_c$; things are boring until failure occurs.  (The most-common numerical symptom is overflow of $e^u$.)  Equation \eqref{liouvillebratu} is a very simple model for combustion of a chemical mixture, and this instability corresponds to a chemical explosion \cite{FrankKameneckij1955}.  However, finding $\lambda_c$ precisely is not easy because \texttt{fas1.py} always initializes at the distant location $w^0=0$.  The behavior of FAS F-cycles is especially nontrivial near the critical $\lambda$ because the critical value is different on coarse grids.  (And apparently sometimes smaller!)  A better strategy for solutions near the critical value, and for parameter studies generally, is ``continuation''.  For example, one might use a saved fine-mesh solution as the initial value in a run with a slightly-different $\lambda$ value.  The new run would only use V-cycles.
\item Equation \eqref{liouvillebratu} is a ``semilinear'' ODE because its nonlinearity occurs in the zeroth-derivative term \cite{Evans2010}.  One might instead solve a ``quasilinear'' equation where the nonlinearity is in the coefficient to the highest-order derivative.  For example, one might try a $p$-Laplacian \cite{Evans2010} extension to the Liouville-Bratu equation:
\begin{equation}
  -\left(|u'|^{p-2} u'\right)' - \lambda e^u = g.  \label{pbratu}
\end{equation}
This equation is the same as \eqref{liouvillebratu} when $p=2$, but for other values $p$ in $(1,\infty)$ the solution is less well-behaved because the coefficient of $u''$ can degenerate or explode.  (A literature at least exists for the corresponding Poisson-like problem with $\lambda=0$ \cite{BarrettLiu1993,Bueler2021}.)  A basic technique is to regularize the leading coefficient with a numerical parameter $\eps>0$: replace $|u'|^{p-2}$ with $\left(|u'|^2+\eps\right)^{(p-2)/2}$.  With such a change, continuation (item \textbf{I}) is also recommended.
\item The most significant extension of \texttt{fas1.py} would be to ``merely'' change from 1D to 2D or 3D.  That is, to change from solving ODEs to solving elliptic PDEs like $-\grad^2 u - \lambda e^u=g$, where $\grad^2$ is the Laplacian operator.  However, doing this in the style of \texttt{fas1.py}, using only NumPy vectors for infrastructure, is not recommended.  Instead, it would be wise to apply an FE library like Firedrake \cite{Rathgeberetal2016} or Fenics \cite{Loggetal2012}, on top of an advanced solver library like PETSc \cite{Balayetal2021,Bueler2021}.  Such libraries imply a substantial learning curve, and their support for FAS multigrid methods is incomplete, but they allow experimentation with higher-order FE spaces and many other benefits.
\end{enumerate}

\section{Conclusion}  \label{sec:conclusion}

Regarding the performance of the tested solvers, we summarize as follows:

\begin{quotation}
\emph{On any mesh of $m$ elements, problem \eqref{feweakform} can be solved nearly to the discretization error of our piecewise-linear FE method by using a single FAS F-cycle, or a few FAS V-cycles.  The work of these methods is $O(m)$ with a small constant; they are optimal.  The faster F-cycle gives textbook multigrid efficiency.  These facts holds for all $m$ up until rounding errors overwhelm the discretization at around $m=10^6$.  By contrast, single-level NGS requires rapidly-increasing numbers of sweeps; its work scales as $O(m^q)$ for $q\gg 1$.}
\end{quotation}

\small

\bigskip
\bibliography{fas}
\bibliographystyle{siam}

\end{document}